\documentclass[11pt]{amsart}
\usepackage{graphicx}

\usepackage{amscd}
\usepackage{mathrsfs}

\theoremstyle{definition}

\def\bkR{\mathbb R}
\def\bkC{{\mathbb C}}
\def\bkN{{\rm \kern.50em \vrule width.05em height1.4ex depth-.05ex \kern-.26em N}}

\def\D{{\mathcal D}}

\newsymbol\subsetneq 2328

\begin{document}
\title[A multiplicative product of distributions and applications]{A multiplicative product of distributions and a class of ordinary differential equations with distributional coefficients}

\author{Nuno Costa Dias and Jo\~ao Nuno Prata}

\address{Departamento de Matem\'atica, Universidade Lus\'ofona de Humanidades e
Tecnologias, Av. Campo Grande 376, 1749-024 Lisboa,
Portugal and Grupo de F\'{\i}sica Matem\'atica,
Universidade de Lisboa, Av. Prof. Gama Pinto 2, 1649-003,
Lisboa, Portugal}
\email{ncdias@meo.pt, joao.prata@mail.telepac.pt}



{\maketitle }

\begin{abstract}
We construct a generalization of the multiplicative product of distributions presented by L. H\"ormander in [L. H\"ormander, {\it The analysis
of linear partial differential operators I} (Springer-Verlag, 1983)]. The new product is defined in the vector space ${\mathcal A}(\bkR)$ of
piecewise smooth functions $f: \bkR \to \bkC$ and all their (distributional) derivatives. It is associative, satisfies the Leibniz rule and
reproduces the usual pointwise product of functions for regular distributions in ${\mathcal A}(\bkR)$. Endowed with this product, the space
${\mathcal A}(\bkR)$ becomes a differential associative algebra of generalized functions. By working in the new ${\mathcal A}(\bkR)$-setting we
determine a method for transforming an ordinary linear differential equation with general solution $\psi$ into another, ordinary linear
differential equation, with general solution $\chi_{\Omega} \psi$, where $\chi_{\Omega}$ is the characteristic function of some prescribed
interval $\Omega \subset \bkR$.
\end{abstract}

\section{Introduction}
Two of the most interesting properties of the space of Schwartz distributions ${\mathcal D}'(\bkR^n)$ over $\bkR^n$ are that the space of
continuous functions $C^{0}(\bkR^n)$ is a subspace of ${\mathcal D}'(\bkR^n)$ and that differentiation is an internal operation in ${\mathcal
D}'(\bkR^n)$ \cite{Schwartz1,Zemanian,Hormander}. On the other hand, its major limitation is that it only displays the structure of a vector
space and not that of an algebra. This has been known since 1954 when L. Schwartz proved that there is {\it no} associative and commutative
algebra of generalized functions
$({\mathcal A}(\bkR^n),+,\circ)$ satisfying the three following properties:\\
\\
(1) The space of distributions ${\mathcal D}'(\bkR^n)$ over $\bkR^n$ is linearly embedded into ${\mathcal A}(\bkR^n)$ and $f(x)\equiv 1$ is the identity in ${\mathcal A}(\bkR^n)$.\\
\\
(2)  The restriction of the product $\circ$ to the set of continuous functions $C^0(\bkR^n)$ reproduces the pointwise product of functions.\\
\\
(3) There exist linear derivative operators $\frac{\partial}{\partial {x_i}}:{\mathcal A}(\bkR^n)\to {\mathcal A}(\bkR^n)$, $(i=1,..,n)$ that:
(a) Satisfy the Leibniz rule and
(b) Their restrictions to ${\mathcal D}'(\bkR^n)$ coincide with the usual distributional derivatives.\\

This theorem became known as the Schwartz impossibility result \cite{Schwartz2}. In spite of it, certain products of distributions emerge rather naturally
in several fields of research in both mathematics and physics. Common examples can be found in quantum electrodynamics and particle physics when dealing with self interacting terms; in hydrodynamics in the formulation of the dynamics of shock waves; and, more generally, whenever one wants
to consider differential equations where either the coefficients or the prospective solutions are non-smooth.

These issues motivated several proposals towards a precise construction of a multiplicative product of distributions. Among these, the most
famous example is due to J. Colombeau \cite{Colombeau0,Colombeau,Colombeau2}, who constructed differential algebras of generalized functions
${\mathcal G}(\bkR^n)$ satisfying conditions (1) and (3) with (2) holding only in $C^{\infty}(\bkR^n)$ (another proposal of a globally defined
product of distributions can be found in \cite{Sarrico1,Sarrico2}). Colombeau algebras are a fine solution for the problem of constructing
associative and commutative algebras ${\mathcal G}(\bkR^n)$ in the situation:
\begin{equation}
C^{\infty}(\bkR^n) \subset C^0(\bkR^n) \subset {\mathcal D}'(\bkR^n)  \subset {\mathcal G}(\bkR^n)
\end{equation}
and satisfying the maximal number of properties stated in the Schwartz impossibility result. Its main disadvantage might be that, in practical
applications, one is required to leave the simpler distributional framework and work in the more involved setting ${\mathcal G}(\bkR^n)$. This
is because, in general, the Colombeau product of two distributions is no longer a distribution. Moreover, $C^0(\bkR^n)$ is not a subalgebra of
${\mathcal G}(\bkR^n)$.

In this paper we want to consider the following problem: Since it is not possible to construct a superspace of
${\mathcal D}'(\bkR^n)$ satisfying the properties (1) to (3) stated in the Schwartz impossibility result, one could, at least, try to construct (the largest possible) subspace of ${\mathcal D}'(\bkR^n)$ that satisfies properties (2) and (3). More precisely:\\

{\bf Problem 1}\\
Determine associative (possible non-commutative) algebras $({\mathcal A}(\bkR^n),+,\star)$ satisfying
\begin{equation}
C^{\infty}(\bkR^n) \subset {\mathcal A}(\bkR^n) \subseteq {\mathcal D}'(\bkR^n)
\end{equation}
and such that:\\
\\
(i) The restriction of the product $\star$ to the set ${\mathcal A}(\bkR^n)\cap C^0(\bkR^n)$ coincides with the usual pointwise product of functions, i.e. $f\star g= f  g$, for all $f,g \in {\mathcal A}(\bkR^n)\cap C^0(\bkR^n)$.\\
\\
(ii) There exist derivative operators in ${\mathcal A}(\bkR^n)$, which are of the form $\frac{\partial}{\partial {x_i}} :{\mathcal A}(\bkR^n)
\to {\mathcal A}(\bkR^n)$ and:

(ii-1) coincide with the restrictions to ${\mathcal A}(\bkR^n)$ of the usual distributional derivatives in ${\mathcal D}'(\bkR^n)$,

(ii-2) satisfy the Leibniz rule.
\\

In this paper we will address this problem in one dimension. We will take ${\mathcal A}(\bkR)$ to be the space $C^{\infty}_p(\bkR)$ of piecewise
smooth functions together with their (distributional) derivatives to all orders, and define a new product $\star$ such that $({\mathcal
A}(\bkR),+,\star)$ becomes an associative (but noncommutative) algebra satisfying properties (i) and (ii) above. The new product $\star$ is a
generalization of the product of distributions with non-intersecting singular supports that was proposed by L. H\"ormander in [pag. 55 of
ref.\cite{Hormander}]. Our task here will consist of extending the domain of H\"ormander's product to include the case of distributions with
intersecting singular supports. Clearly, this is not possible for the entire set ${\mathcal D}'(\bkR)$, but we will prove that it is possible in
${\mathcal A}(\bkR)$.

It is important to remark that:\\
(a) ${\mathcal A}(\bkR)$ is not a subalgebra of ${\mathcal G}(\bkR)$ because the product $\star$ is not the restriction of the Colombeau product $\circ$ to ${\mathcal A}(\bkR)$.  In fact, Colombeau's product is not an inner operation in ${\mathcal A}(\bkR)$ and it does not even satisfy $f\circ g \in {\mathcal D}'(\bkR)$ for all $f,g \in {\mathcal A}(\bkR)$.\\
(b) In many specific problems, where products of distributions are present, what is at stake are products involving piecewise
$C^{\infty}(\bkR)$-functions and  distributional derivatives of these functions. This is the typical case when considering differential
equations with distributional coefficients or non-linear terms.
For these cases, the algebra ${\mathcal A}(\bkR)$ displays almost optimal properties. It provides a sufficiently general setting without ever leaving the space ${\mathcal D}'(\bkR)$.\\
(c) The main limitation of ${\mathcal A}(\bkR)$ is that it is defined over $\bkR$, only. Moreover, the details of the derivation of the product $\star$ indicate that the extension of its domain to functionals of several variables might not be easy. We hope to study this issue in a forthcoming paper.\\

In the second part of this paper, we use our product of distributions to address the following problem:\\

{\bf Problem 2}\\
Let $\Omega \subset \bkR$ be an open interval and let $\chi_{\Omega}$ be the characteristic function of $\Omega$ (i.e. $\chi_{\Omega}(x)=1$ for $x\in \Omega$ and zero otherwise).\\
Let us consider a generic, linear and ordinary differential equation and let $\psi \in C^{\infty}(\bkR)$ be its general solution on $\bkR$.\\
The problem we want to address is that of deriving a new, also linear and ordinary differential equation, but displaying the general solution $\chi_{\Omega}\psi$.\\

The crux of the matter here is that it is the {\it general solution} (i.e. irrespectively of the boundary conditions) that has to be confined to
the interval $\Omega$. The case where we only impose the confinement of a particular solution can be solved by adding to the original
differential equation a suitable non-homogenous distributional term. Differential equations of this type and their relation with the Green's
function theory and with functions with jump discontinuities have been investigated in \cite{Kanwal,Estrada,Pan,Stakgold}. In general the
distributional term produces a jump discontinuity of fixed strength that is able to confine a particular, but not all solutions \cite{Kanwal}.
In several cases this is an important limitation. For instance, if we want to construct a (version of a given) differential operator whose
eigenfunctions are all confined to an interval $\Omega$ (a relevant problem, for example, for the global formulation of quantum systems with
boundaries \cite{Isham}) the "distributional non-homogeneous term" method is useless.

By working in ${\mathcal A}(\bkR)$ we will provide a general method for constructing the differential equations, solutions of Problem 2. Their
most interesting feature is that they display an extra term written in terms of a distributional coefficient. This is why they require the use
of a product of distributions. The new differential equations yield a global formulation for a class of differential problems defined on bounded
domains and may find interesting applications in those fields of mathematical physics and dynamical systems where global analysis is required
\cite{Isham}. An example of such applications is the derivation of a {\it globally defined} time independent Schr\"odinger equation for systems
with boundaries \cite{Isham,Dias1,Piotr,Dias2} and the related issue of quantum confinement in several non-local formulations of quantum
mechanics, such as the deformation and the De Broglie-Bohm formulations \cite{Walton,Dias3,Dias4}. These problems were studied in
\cite{Dias1,Dias2,Dias3} using an approach related to the results that will be presented here.

This paper is organized as follows: In section 2 we introduce the notation, provide a concise review of some relevant results on Schwartz's
distributions and define the space ${\mathcal A}(\bkR)$. In section 3 we define a product in ${\mathcal A}(\bkR)$ and study its main properties.
In section 4 we address the problem of constructing globally defined differential equations displaying solutions of confined support.

\section{Schwartz distributions and the space ${\mathcal A}(\bkR)$}

In this section we introduce the notation, review some classical results in the theory of distributions and define the space ${\mathcal
A}(\bkR)$ (in the next section we will prove that ${\mathcal A}(\bkR)$ is an algebra). Let $\Omega \subseteq \bkR$ be an open interval,
$C^m(\Omega)$ the vector space of complex functions on $\Omega$ which are $m$-times continuously differentiable, and let $C^m_p(\Omega)$ be the
vector space of piecewise $C^m$-functions on $\Omega$, i.e. of functions that are $m$-times continuously differentiable, except on a finite set,
where they and all their derivatives up to order $m$ have finite left and right limits. Let also $\D(\Omega )$ be the vector space of infinitely
smooth, complex valued functions with support on a compact subset of $\Omega$. The set $\D(\Omega)$ is usually named the set of test functions.
The set of Schwartz distributions $\D'(\Omega )$ is the topological dual of $\D(\Omega)$. Its elements are linear and continuous functionals $F$
acting on $\D(\Omega)$ \cite{Zemanian}:
\begin{equation}
F:\D(\Omega )  \to {\bkC}; \quad t \to F(t) \equiv <F,t>
\end{equation}
where $F(t)$ and $<F,t>$ are two alternative notations for the action of $F$ on $t$.

In our notation the lowercase roman letters from the middle of the alphabet ($f,g,h...$) will be used for arbitrary functions, those from the
end of the alphabet ($t,s,r,...$) will be reserved for test functions and the upper case roman letters from the middle of the alphabet
($F,G,J,K$) will be used for distributions. $H$ is the Heaviside step function.
By $\frac{d^n}{dx^n} f$ and  $f^{(n)}$ we denote the $n$th-order distributional derivative of $f$. The notation $f',f''$ will also be used for the first and second order derivatives.

A distribution $F$ is regular if it acts as $<F,t>=\int \, f(x)t(x)dx$ for some locally integrable function $f$.
If $f$ is continuous we shall make the identification $F=f$ because $f$ is determined by $F$ uniquely.
Furthermore, the restriction of a distribution $F:\D (\Omega)\to \bkC$ to an open set $\Xi \subset \Omega$ is the distribution $F_{\Xi}:\, \D(\Xi) \to \bkC ; \, t\to <F_{\Xi},t>=<F,t>$.  Two distributions $F,G \in \D'(\Omega )$ are identical on some open set $\Xi \subset \Omega$ iff $F_{\Xi}=G_{\Xi}$. The null set of $F$ (denoted null $F$) is the largest open set ${\Xi }$ where $F_{\Xi }=0$. The complement of this set is the support of $F$, denoted supp $F$. The singular support of $F$ is defined as the complement of the largest open set ${\Xi }$ where $F=f$ for some $f\in C^{\infty}({\Xi })$. The standard notation for this set is sing supp $F$.
Another important notion is that of (weak) convergence in $\D'(\Omega)$. A sequence of distributions $(F_n)_{1\le n <\infty }$ is said to converge (weakly) to $F\in \D'(\Omega)$ iff $\lim_{n \to \infty} <F_n,t>=<F,t>$ for all $t\in \D(\Omega)$.

The following definition is very important for the sequel:\\
\\
{\bf Definition 2.1: The space ${\mathcal A}(\Omega)$}\\
{\it The space of special distributions ${\mathcal A}(\Omega)$ on $\Omega$ is the space of piecewise smooth functions $C^{\infty}_p(\Omega)$ (regarded as distributions) together with their distributional derivatives to all orders.}\\
\\
{\bf Remark 2.2.} All functions $f \in C_p^{\infty}(\Omega)$ can be associated with a distribution through the prescription $<f,t>= \int \, f(x)t(x)dx$. In this sense, $C^{\infty}_p(\Omega)\subset \D'(\Omega)$ and so ${\mathcal A}(\Omega)$ is a subset of $\D'(\Omega)$. One easily realizes that ${\mathcal A}(\Omega)$ is, in fact, a linear subspace of $\D'(\Omega)$. It also follows from its definition that ${\mathcal A}(\Omega)$ is closed under differentiation. In the next section we will define a product $\star$ in ${\mathcal A}(\bkR)$ and prove that $({\mathcal A}(\bkR),+,\star)$ is an associative, differential algebra.\\
\\
{\bf Remark 2.3.} For each $f\in C_p^{\infty}(\Omega)$ there is a finite set $V_f=\{x_1 < x_2 <...<x_N\} \subset \Omega$ such that $f \in C^{\infty}(\Omega \backslash V_f)$. Let $I_i=]x_i,x_{i+1}[$, $i=1,..,N-1$; $I_0=]-\infty,x_1[ \cap \Omega$ and $I_N=]x_N, +\infty[ \cap \Omega$. Since $f$ and all its derivatives display finite limits at $x_i$, $i=1,..,N$ there is (by Whitney's extension theorem  \cite{Whitney,Fleming}) a collection of functions $f_i \in C^{\infty}(\Omega)$, $i=0,..,N$ such that their restrictions to $I_i$ coincide with those of $f$, i.e:
\begin{equation}
f(x)= f_i(x)\quad \mbox{if}\quad x\in I_i \quad ;\quad i=0,...,N
\end{equation}
Hence, on $\Omega \backslash V_f$ we have:
\begin{equation}
f= \sum_{i=0}^N \chi_i f_i
\end{equation}
where $\chi_i(x)=1$ if $x\in I_i$ and $\chi_i(x)=0$ if $x\notin I_i$ is the characteristic function of $I_i$. Conversely, if $f$ is of the form (5) then $f \in C_p^{\infty}(\Omega)$. These statements are also valid in distributional sense, i.e.
$f\in C_p^{\infty}(\Omega)$ (regarded as a subset of ${\mathcal D}'(\Omega)$ - Remark 2.2) iff there is a finite set $V_f \subset \Omega$ and a collection of functions $f_i \in C^{\infty}(\Omega)$, $i=0,..,N=\sharp V_f$ such that, in the sense of distributions, $f$ can be written in the form (5).\\
\\
{\bf Remark 2.4.} It is clear that $C^{\infty}(\Omega) $ is a subspace of ${\mathcal A}(\Omega)$, but $C^0(\Omega)$ is not. For instance
$f(x)=\sqrt{|x|} \notin {\mathcal A}(\bkR)$. The Dirac delta distribution $\delta$ and all its derivatives are elements of ${\mathcal
A}(\Omega)$.
\\

Lastly, we review the content of two well-known theorems, which will be important in the sequel. Proofs may be found in \cite{Zemanian}.
The first one concerns the (re)construction of a global distribution from a set of local ones:\\
\\
{\bf Lemma 2.5.} {\it Let $N\in \bkN$ and let $\Omega_k \subset \Omega ,\,  k=1,..,N$ be a finite collection of open sets covering $\Omega$. Let $F_k \in \D'(\Omega_k)$, $k=1,..,N$ be a collection of distributions such that $F_i=F_j$ in $\Omega_i\cap \Omega_j$ for all $i,j=1,..,N$. Then there is one and only one $F\in \D'(\Omega)$ such that $F_{\Omega_k}=F_k$ for all $k=1,..,N$.}\\

The second one is a consequence of the completeness of $\D'(\Omega)$ with respect to the weak convergence:\\
\\
{\bf Lemma 2.6.} {\it Let $F_{\epsilon}\in \D'(\Omega)$ be a one parameter family of distributions with $\epsilon >0 $. If $\lim_{\epsilon \to 0^+} <F_{\epsilon},t>$ exists for every $t\in \D(\Omega)$ then:
\begin{equation}
F:\D(\Omega) \longrightarrow \bkC,\, t \longrightarrow <F,t>\equiv \lim_{\epsilon \to 0^+}<F_{\epsilon},t>
\end{equation}
is an element of $\D'(\Omega)$.}

\section{The algebra of special Distributions and multiplicative products}

In this section we will provide two alternative characterizations of the space ${\mathcal A}(\bkR)$ (Theorems 3.1 and 3.2), study H\"ormander's
product of distributions with non intersecting singular supports (Definition 3.4 to Theorem 3.7) and make the proposal of the new product
$\star$ (Definition 3.10 and Theorems 3.13 to 3.16). These results prove that $({\mathcal A}(\bkR),+,\star )$ is indeed a differential
associative algebra and a solution of {\it Problem 1}. To make the presentation simpler, we will assume that $\Omega = \bkR$, but our results
are still valid for an arbitrary open interval $\Omega \subset \bkR$.

The next theorem follows from Remark 2.3 and provides an important characterization of ${\mathcal A}(\bkR)$:\\
\\
{\bf Theorem 3.1.} {\it $F\in {\mathcal A}(\bkR)$ iff there is a finite set of real numbers $V_F\subset \bkR$ and two Schwartz distributions $f$
and $\Delta^{(F)}$, such that:
\begin{equation}
F=f+\Delta^{(F)}
\end{equation}
where $f\in C_p^{\infty}(\bkR) \cap C^{\infty}(\bkR \backslash V_F) $, and $\Delta^{(F)}$ is of the form:
\begin{equation}
\Delta^{(F)}=\sum_{w\in V_F} \Delta^{(F)}_w =
\sum_{w\in V_F} \sum_{k=0}^N C_k(w) \delta^{(k)}(x-w)
\end{equation}
for some $N\in \bkN_0$ and $C_k(w) \in \bkC$.}
\\
\\
{\bf Proof.} If $F\in {\mathcal A}(\bkR)$, then $F$ is a $n$th-order (distributional) derivative of some function $g \in C^{\infty}_p(\bkR)$
(cf. Definition 2.1). It follows from Remark 2.3 that, in the sense of distributions, $g$ can be written as:
\begin{equation}
g=g_0 + \sum_{w\in V_g} H(x-w) g_w
\end{equation}
where $H$ is the Heaviside step function, $g_0,g_w \in C^{\infty}(\bkR)$ and $V_g$ is a finite set of real numbers. We then have:
\begin{equation}
F = \frac{d^n}{dx^n} g = g_0^{(n)} + \sum _{w \in V_g} \sum_{p=0}^n \left(^n_p \right) H^{(n-p)}(x-w) g_w^{(p)}
\end{equation}
where $(^n_p)= \frac{n!}{(n-p)!p!}$ is the binomial coefficient. It is easy to realize that $F$ can be cast in the form (7) by setting:
\begin{eqnarray}
 f  & = & g_0^{(n)} + \sum _{w \in V_g} H(x-w) g_w^{(n)}   \\
 \Delta^{(F)} & = &  \sum_{w\in V_g}  \sum_{p=0}^{n-1} \left(^n_p \right) g_w^{(p)}(x) \delta^{(n-1-p)}(x-w) \nonumber \\
 & = &  \sum_{w \in V_g}  \sum_{p=0}^{n-1} g_w^{(p)} (w) \delta^{(n-1-p)}(x-w)  \nonumber
\end{eqnarray}
and it is clear that i) $f\in C_p^{\infty}(\bkR) \cap C^{\infty}(\bkR\backslash V_g)$ and that ii) $\Delta^{(F)}$ is of the form (8)).\\
Conversely, if $F$ is given by eqs.(7,8) then $F$ is the $(N+1)$th-order derivative of:
\begin{equation}
g=g_0+ \sum _{w \in V_F} \sum_{k=0}^N \frac{C_k(w)}{(N-k)!} H(x-w) (x-w)^{N-k}
\end{equation}
where:
\begin{equation}
g_0(x) = \left\{ \begin{array}{l}
\int_a^{x} dx_0 \, f(x_0) \quad , \quad N = 0 \\
\\
\int_a^{x} dx_0 \int_a^{x_0} dx_1 \int_a^{x_1} dx_2 .... \int_a^{x_{N-1}} dx_N \, f(x_N) \quad , \quad N \ge 1
\end{array} \right.
\end{equation}
for some $a \in \bkR$. Since $g \in C_p^{\infty}(\bkR)$ we have $F \in {\mathcal A}(\bkR)$ which concludes the proof.$_{\Box}$
\\

The next theorem characterizes ${\mathcal A}(\bkR)$ as the set of Schwartz distributions with finite singular support which, in addition, are identical to a $C^{\infty}(\bkR)$ function on every open interval not intersecting its singular support. The theorem is not essential for the sequel and can be skipped.\\
\\
{\bf Theorem 3.2.} {\it A distribution $F\in \D'(\bkR)$ belongs to ${\mathcal A}(\bkR)$ iff:
\\
(a) The singular support of $F$ is a finite set, and
\\
(b) If ${\Xi }\subset \bkR$ is an open interval satisfying ${\Xi } \cap \,$}sing supp $F = {\emptyset } ${\it, then there is a function
$g\in C^{\infty}(\bkR)$ such that the identity "$F=g$ on ${\Xi }$" is valid in the distributional sense.}
\\
\\
{\bf Proof.} We will prove this result by using the definition of ${\mathcal A}(\bkR)$ that follows from Theorem 3.1. If $F\in \D'(\bkR)$
satisfies the condition (a) above, we may set $V_F=$sing supp $F$ which is then a finite set. Moreover, the restriction $F_{\bkR \backslash
V_F}=f_{\bkR \backslash V_F}$ for some $f \in C^{\infty}(\bkR \backslash V_F)$. We then have $(F-f)_{\bkR \backslash V_F}=0$ and so supp$(F-f)
\subseteq V_F$. Let $\Delta^{(F)}=F-f$. Since $\Delta^{(F)}$ is supported on the finite set $V_F$ then (by a well-known textbook result
\cite{Zemanian}) $\Delta^{(F)}$ is a finite linear combination of Dirac deltas and their derivatives, i.e.:
\begin{equation}
\Delta^{(F)}=\sum_{w\in V_F} \Delta^{(F)}_w  \quad \mbox{where} \quad \Delta^{(F)}_w(x)=\sum_{k=0}^N C_k(w) \delta^{(k)}(x-w), \quad N\in
\bkN_0,\, C_k(w) \in \bkC
\end{equation}
in agreement with eq.(8). Hence, $F=f+\Delta^{(F)}$ with $\Delta^{(F)}$ given by (8) and $f \in C^{\infty}(\bkR \backslash V_F)$. To prove that
$F\in {\mathcal A}(\bkR)$, by Theorem 3.1, we still have to prove that $f \in C_p^{\infty}(\bkR)$. Let ${\Xi }=]x_0,y_0[$ where $x_0,y_0$ are
two arbitrary consecutive elements of $V_F$. Then $\Xi \cap V_F=\emptyset $ and from condition (b) above $F_{\Xi }=g_{\Xi }$ for some $g \in
C^{\infty}(\bkR )$. Since it is also true that $F_{\Xi }=f_{\Xi }$ it follows that $f_{\Xi }=g_{\Xi }$ and so $\lim_{x \to x_0^+}f^{(k)}=\lim_{x
\to x_0^+}g^{(k)}$ and $\lim_{x \to y_0^-}f^{(k)}=\lim_{x \to y_0^-}g^{(k)}$ exist and are finite for all $k\in \bkN _0$. This proves that $f\in
C_p^{\infty}(\bkR)$ and so, by Theorem 3.1, that $F \in {\mathcal A}(\bkR)$.
\\
Conversely, if $F\in {\mathcal A}(\bkR)$ then $F$ is given by eqs.(7,8) for some finite set $V_F$ and some $f\in C_p^{\infty}(\bkR)$. Hence,
sing supp $F$ is a finite set. Moreover, on an arbitrary open interval ${\Xi }\subset \bkR$ satisfying ${\Xi } \cap \,$sing supp $F= {\emptyset
}$, we have $F_{\Xi }=f_{\Xi }$. Since $f$ is infinitely smooth on ${\Xi }$ and $f^{(k)}$ displays finite lateral limits at the boundaries of
${\Xi }$ (for all $k\in \bkN_0$), there is (by Whitney's extension theorem \cite{Whitney}) a function $g\in C^{\infty}(\bkR)$ such that $f_{\Xi
}=g_{\Xi }$ and so also  $F_{\Xi }=g_{\Xi }$.
This proves that $F$ also satisfies condition (b) above.$_{\Box}$\\

Before proceeding, let us make the following remark concerning the relation between $V_F$ and sing supp $F$:\\
\\
{\bf Remark 3.3.} We will always assume that $V_F$ is a finite superset of sing supp $F$ (but not necessarily identical to sing supp $F$). This is not contradictory with Theorems 3.1 and 3.2 (since sing supp $F$ is always a finite set and we may always set $\Delta^{(F)}_w=0$ for some $w \in V_F$). The reason why we do not simply make $V_F=$sing supp $F$ (although, most of the time, we assume this to be the case) is that the extra freedom will allow us to simplify the notation when deriving the explicit expression for the new product of distributions. \\

The aim now is to introduce a multiplicative product in the space ${\mathcal A}(\bkR)$. We start by considering the (restriction to ${\mathcal A}(\bkR)$ of the) product of distributions with {\it non intersecting singular supports} that was proposed by L. H\"ormander in [pag. 55, \cite{Hormander}].\\
\\
{\bf Definition 3.4: The H\"ormander product $\cdot$ }
\\
{\it Let $F,G \in {\mathcal A}(\bkR)$ be such that $V_F \cap V_G= {\emptyset } $. Let $\{\Omega_w \subset \bkR, w\in V_F \cup V_G\}$ be a finite
covering of $\bkR$ by open sets satisfying $w'\notin \Omega_w$, $\forall w\not=w' \in V_F \cup V_G$. Let us also introduce the compact notation
$F_w=F_{\Omega_w}$, $G_w=G_{\Omega_w}$ to designate the restrictions of $F$ and $G$ to the set $\Omega_w$. Then $F\cdot G$ is defined by:
\begin{equation}
F\cdot G: \quad (F\cdot G)_{\Omega_w}=F_wG_w
\end{equation}
where $F_wG_w$ denotes the usual product of a distribution by an infinitely smooth function.}\\

The Remarks 3.5 and 3.6 will show that the Definition 3.4 is consistent. They are included for completion. Remark 3.5 also provides the explicit form of $F_wG_w$. This result will be used in Theorem 3.7 to obtain an explicit expression for $F\cdot G$.
\\
\\
{\bf Remark 3.5.} The product $F_wG_w$ in eq.(15) is well defined. Indeed, let us write $F,G\in {\mathcal A}(\bkR)$ in the form given by eq.(7),
i.e. $F=f+\Delta^{(F)}$ and $G=g+\Delta^{(G)}$. For each $w \in V_F \cup V_G $ either $F_w$ or $G_w$ is a regular distribution associated to a
function in $C^{\infty}(\Omega_w)$:
\begin{eqnarray}
 w\in V_F &\Longrightarrow &\left\{
\begin{array}{l}
F_w =  \Delta^{(F)}_w+f_w ; \\
f_w \in C^{\infty}(\Omega_w \backslash \{ w\} ) \cap C^p (\Omega_w)\\
\\
G_w =  g_w \in C^{\infty}(\Omega_w)
\end{array} \right. \qquad \mbox{and} \\
w\in V_G &\Longrightarrow & \left\{
\begin{array}{l}
F_w  =  f_w \in C^{\infty}(\Omega_w)\\
 \\
 G_w =  \Delta^{(G)}_w+g_w ;\\
 g_w \in C^{\infty}(\Omega_w \backslash \{ w\} ) \cap C^p (\Omega_w)
\end{array} \right. \nonumber
\end{eqnarray}
where $g_w=g_{\Omega_w}$ and $f_w=f_{\Omega_w}$ are the restrictions of $g$ and $f$ to the open set ${\Omega_w}$. It follows that
the product $F_wG_w$ is well defined in the usual sense of a product of a distribution by an infinitely smooth function. We have:
\begin{equation}
F_wG_w= \left\{
\begin{array}{l}
g_w\Delta^{(F)}_w+g_wf_w, \quad w\in V_F \\
\\
f_w\Delta^{(G)}_w+g_wf_w, \quad w\in V_G \\
\end{array} \right.
\end{equation}
and so $F_wG_w \in \D'(\Omega_w)$.\\
\\
{\bf Remark 3.6.} The definition of $F\cdot G$ satisfies the prerequisites of Lemma 2.5. Indeed, on an arbitrary open set
${\Xi }_{ww'}=\Omega_w \cap \Omega_{w'}$ (with $w\not=w'$) we have $(F_wG_w)_{{\Xi }_{ww'}}=(fg)_{{\Xi }_{ww'}}=(F_{w'}G_{w'})_{{\Xi }_{ww'}}$. This is because ${{\Xi }_{ww'}}\cap (V_F\cup V_G)= {\emptyset } $. Lemma 2.5 then implies that $F \cdot G$ is a Schwartz distribution, uniquely defined by eq.(15).\\

We have not yet proved that $F\cdot G$ is independent of the particular covering $\{\Omega_w, \, w \in V_F \cup V_G \}$ used to define it. This result will be a corollary of the next theorem where we obtain an explicit expression for $F \cdot G$:\\
\\
{\bf Theorem 3.7.} {\it Let $F,G\in {\mathcal A}(\bkR)$ (such that $V_F \cap V_G =\emptyset$) be written in the form (7,8), i.e. $F=f + \sum_{w
\in V_F} \Delta^{(F)}_w$ and $G=g + \sum_{w \in V_G} \Delta^{(G)}_w$. Then $F\cdot G$ is given by:
\begin{equation}
F\cdot G= fg + \sum_{w \in V_F} \tilde g_w\Delta^{(F)}_w+\sum_{w \in V_G} \tilde f_w\Delta^{(G)}_w
\end{equation}
where $\tilde g_w$ and $\tilde f_w$ are $C^{\infty}(\bkR)$-extensions of the restrictions $g_w=g_{\Omega_w}$ and $f_w=f_{\Omega_w}$, respectively; and
$\{\Omega_w \subset \bkR, \, w\in V_F \cup V_G \}$ is a finite covering of $\bkR$ satisfying the requisites of Definition 3.4.}
\\
\\
{\bf Remark 3.8.} The extensions $\tilde g_w$ and $\tilde f_w$ are $C^{\infty}(\bkR)$-functions such that, for each $w\in V_F$ we have
$\tilde g_w(x)=g(x)$ if $x\in \Omega_w$; and for each $w\in V_G$, $\tilde f_w(x)=f(x)$ if $x\in \Omega_w$. That these functions exist follows from Whitney's extension theorem \cite{Whitney,Fleming} and the fact that $f,g \in C_p^{\infty}(\bkR)$ and so, for each $w\in V_F$ (respectively $w\in V_G$), the restrictions $g_w \in C^{\infty}(\Omega_w)$ -cf. eq.(16)- (respectively $f_w \in C^{\infty}(\Omega_w)$) and all its derivatives display finite limits at the boundaries of $\Omega_w$.
Hence, the right hand side of eq.(18) is a well defined Schwartz distribution. Let us denote it by $J$.
\\
\\
{\bf Remark 3.9.} We now prove that $J$ is independent of the particular finite covering of $\bkR$
(satisfying the conditions in Definition 3.4) and of the particular extensions of $g_w$ and $f_w$ that were used to define it. Indeed, if
$\{\Omega_w, \, w\in V_F \cup V_G\}$ and $\{\Omega'_w\, w\in V_F \cup V_G\}$ are two such coverings and $\tilde g_w,\tilde f_w$ and $\tilde g'_w,\tilde f'_w$ the associated functions then $w
\notin $ supp$(\tilde f'_w-\tilde f_w)$, $\forall w \in V_G$ and $w \notin $ supp$(\tilde g'_w-\tilde g_w)$, $\forall w\in V_F$. Hence:
$$
\mbox{supp} \, (\tilde g'_w-\tilde g_w) \cap \,\mbox{supp} \, \Delta^{(F)}_w= {\emptyset } , \quad \forall w \in V_F
$$
and
$$
\mbox{supp} \, (\tilde f'_w-\tilde f_w) \cap \,\mbox{supp} \, \Delta^{(G)}_w= {\emptyset } ,  \quad \forall w \in V_G
$$
and so
\begin{equation}
\sum_{w \in V_F} \tilde g_w\Delta^{(F)}_w+\sum_{w \in V_G} \tilde f_w\Delta^{(G)}_w- \left[\sum_{w \in V_F}
\tilde g'_w\Delta^{(F)}_w+\sum_{w \in V_G} \tilde f'_w\Delta^{(G)}_w \right]=0
\end{equation}
Therefore $J$ is a well defined distribution that only depends of $F$ and $G$.
\\
\\
{\bf Proof of Theorem 3.7.}
It is straightforward to realize that the restriction of $J$ to $\Omega_w$ is $F_wG_w$ (given by eq.(17)). Since the restrictions to the open covering $\{\Omega_w,\, w\in V_F\cup V_G\}$ define a unique global distribution (Lemma 2.5) which, by definition, is $F\cdot G$ (cf. eq.(15)) the identity $F\cdot G=J$ holds.$_{\Box}$\\

Two simple corollaries follow directly from the theorem:\\
\\
{\bf Corollary 3.10.} $F \cdot G$ is independent of the particular open covering (Definition 3.4) used to define it.
\\
\\
{\bf Corollary 3.11.} If $F,G \in {\mathcal A}(\bkR)$ (such that $V_F \cap V_G =\emptyset $) then $F\cdot G \in {\mathcal A}(\bkR)$ and sing supp $({F\cdot G}) \subseteq (V_F \cup V_G)$.\\

This concludes the study of H\"ormander's product. The next definition
introduces the new product $\star$, extending H\"ormander's product to the case of distributions with intersecting singular supports.\\
\\
{\bf Definition 3.12: The product $\star$} \\
{\it Let $F,G \in {\mathcal A}(\bkR)$. The multiplicative product $\star$ in ${\mathcal A}(\bkR)$ is defined by:
\begin{equation}
F\star G= \lim_{\epsilon \to 0^+}F \cdot G^{\epsilon}
\end{equation}
where $G^{\epsilon}(x)=G(x+\epsilon)$ is the translation of $G$ by $\epsilon $ and the limit is taken in the weak sense.}
\\

The two following Remarks discuss the consistency of Definition 3.12. Then, in Theorem 3.15 we will prove that $\star$ is an internal operation
in ${\mathcal A}(\bkR)$ and determine the explicit form of $F\star G$ for general $F,G \in {\mathcal A}(\bkR)$.
\\
\\
{\bf Remark 3.13.} Since $V_F$ and $V_G$ are finite sets of real numbers there is always a real number $\sigma >0$ such that $V_{F}\cap V_{G^{\epsilon}}= {\emptyset } $, $\forall \epsilon \in ]0,\sigma[$; notice that $V_{G^{\epsilon}}=\{x-\epsilon, \, x\in V_G \}$ and so we may set, for instance, $\sigma = \mbox{min} \left(\{|y-x|, x \in V_F, y \in V_G \} \backslash \{0\} \right)$ if this minimum exists and $\sigma=1$ if this is not the case (because $V_F=V_G$ and $\sharp V_F=1$). Hence the product $F \cdot G^{\epsilon}$ is well defined for $\epsilon \in ]0,\sigma[$.\\
Notice that, in general, $F\cdot G$ is not well defined since sing supp $F\cap $ sing supp $G\not=\emptyset$. However, when that is the case, we will easily find (from Theorem 3.15) that $F \star G=F\cdot G$.\\
\\
{\bf Remark 3.14.} The functional $F\star G$ acts as:
\begin{equation}
<F\star G,t>=\lim_{\epsilon \to 0^+} <F\cdot G^{\epsilon},t>
\end{equation}
and if its domain is $\D(\bkR)$ - i.e. if the limit on the right hand side exists for all $t\in \D(\bkR)$ - then $F\star G$ is a Schwartz
distribution (cf. Lemma 2.6). The next theorem will prove that this is the case for general $F,G \in {\mathcal A}(\bkR)$ and furthermore that
$F\star G \in {\mathcal A}(\bkR)$.
\\

Before we proceed we need to introduce some notation. Let $F,G \in {\mathcal A}(\bkR)$ and (for notational convenience) let us set
$V_F=V_G=$sing supp $F \, \cup \,$sing supp $G$. Let also $N=\sharp (V_F)$ and $x_k\in V_F$, $k=1,..,N$, $x_{k+1}>x_k$ be the array of elements
of $V_F$ (or $V_G$). Hence, $V_F=V_G=\{x_1,...,x_N\}$ and the distributions $F$ and $G$ can then be written in the form (cf. eqs(7,8)):
\begin{equation}
\left\{ \begin{array}{l}
F= f+\sum_{k=1}^N  \Delta^{(F)}_{x_k}  \\
\\
G=g+\sum_{k=1}^N  \Delta^{(G)}_{x_k}
\end{array}
\right.
\end{equation}
where we set $\Delta^{(F)}_{x_k}=0$ for $x_k\in V_F {\backslash }\, $sing supp $F$ and $\Delta^{(G)}_{x_k}=0$ for $x_k\in V_G {\backslash } \,
$sing supp $G$. Moreover, $f,g \in C_p^{\infty}(\bkR) \cap C^{\infty}(\bkR \backslash V_F)$ and so (cf. Remark 2.3) they can be written in terms of a set of
$C^{\infty}(\bkR)$-functions $f_k,g_k$; $k=0,..,N$:
\begin{eqnarray}
f(x) & = & H(x_1-x)f_0(x) + \sum_{k=1}^{N-1} H(x-x_k) H(x_{k+1}-x) f_k(x) + H(x-x_N)f_N(x)  \nonumber \\
\\
g(x) & = & H(x_1-x)g_0(x) + \sum_{k=1}^{N-1} H(x-x_k) H(x_{k+1}-x) g_k(x) + H(x-x_N)g_N(x).\nonumber
\end{eqnarray}
We can now state the new theorem:\\
\\
{\bf Theorem 3.15.} {\it Let $F,G \in {\mathcal A}(\bkR)$. Then $F\star G \in {\mathcal A}(\bkR)$ and its explicit expression is given by:
\begin{equation}
F\star G= fg + \sum_{k=1}^N \left[ g_k
\Delta^{(F)}_{x_k}  +f_{k-1} \Delta^{(G)}_{x_k} \right]
\end{equation}}
for $F,G$ written in the form (22,23).
\\
\\
{\bf Proof.} The first step is to obtain an explicit expression for $F\cdot G^{\epsilon}$ (cf. Definition 3.12). Let then $F,G$ be written in the form (22,23).
It follows from eq.(22) that $G^{\epsilon}(x) =G(x+\epsilon)$ is given by:
\begin{equation}
G^{\epsilon}=g^{\epsilon} + \sum_{k=1}^N  \Delta^{(G)}_{x_k-\epsilon}
\end{equation}
where $g^{\epsilon}(x)=g(x+\epsilon)$ and $\Delta^{(G)}_{x_k-\epsilon}(x)=\Delta^{(G)}_{x_k}(x+\epsilon)$ is a singular distribution with
support on $\{x_k-\epsilon\}$. We also have $V_F=V_G=\{x_k, k=1,..,N\}$ and so $V_{G^{\epsilon}}=\{x_k-\epsilon, \, k=1,..,N\}$.
Since there is a $\sigma >0$ such that $V_F \cap V_{G^{\epsilon}}= {\emptyset } $ for $0<\epsilon <\sigma $ (cf. Remark 3.13), we can use eq.(18)  to calculate $F\cdot G^{\epsilon}$ explicitly for all $\epsilon \in ]0,\sigma[$:
\begin{eqnarray}
F\cdot G^{\epsilon} &= & fg^{\epsilon} + \sum_{w \in V_F} \tilde{(g^{\epsilon})}_w \Delta^{(F)}_w+\sum_{w \in V_{G^{\epsilon}}} \tilde
f_w\Delta^{(G)}_w \\
&= & fg^{\epsilon} + \sum_{k=1}^N \tilde{(g^{\epsilon})}_{x_k} \Delta^{(F)}_{x_k}+\sum_{k=1}^N \tilde
f_{(x_k-\epsilon)}\Delta^{(G)}_{(x_k-\epsilon)} \nonumber
\end{eqnarray}
where $\tilde{(g^{\epsilon})}_{x_k}$ and $\tilde f_{(x_k-\epsilon)}$ are the $C^{\infty}(\bkR)$-extensions of the restrictions ${(g^{\epsilon})}_{x_k}={(g^{\epsilon})}_{\Omega_{x_k}}$ and $ f_{(x_k-\epsilon)}=f_{\Omega_{(x_k-\epsilon)}}$, respectively. As usual the open sets $\Omega_{x_k}$ and $\Omega_{(x_k-\epsilon)}$ belong to a finite covering $\{\Omega_w \subset \bkR; \, w \in V_F \cup V_{G^{\epsilon}} \}$ of $\bkR$, satisfying the conditions stated in Definition 3.4. That is: i) $\cup_{k=1}^N (\Omega_{x_k} \cup \Omega_{(x_k-\epsilon)})=\bkR$ and ii) $w' \notin \Omega_w $, $\forall w\not=w' \in (V_F\cup V_{G^{\epsilon}})$. There are, of course, many possible coverings of this sort but the expression (26) is independent of the particular one used to define it (cf. Remark 3.9). One possibility is (let $x_0=-\infty$ and $x_{N+1}=+\infty$):
\begin{equation}
\Omega_{x_k} = ]x_k-\epsilon,x_{k+1}-\epsilon[ \quad \mbox{and} \quad
\Omega_{x_k-\epsilon} = ]x_{k-1},x_{k}[,\quad ,k=1,..,N
\end{equation}
It then follows from eq.(23) that:
\begin{equation}
(g^{\epsilon})_{\Omega_{x_k}}=(g^{\epsilon}_k)_{\Omega_{x_k}} \quad \mbox{and} \quad
(f)_{\Omega_{(x_k-\epsilon)}}=(f_{k-1})_{\Omega_{(x_k-\epsilon)}}\quad , \quad k=1,..,N
\end{equation}
where $g^{\epsilon}_k$ is defined by
$g^{\epsilon}_k(x)=g_k(x+\epsilon)$ and $g_k,\, f_k$ were introduced in eq.(23). Since $g_k$ and $f_k$ belong to $C^{\infty}(\bkR)$ for all $k=0,..,N$ we may set (cf. eq.(26)):
$$
\tilde{(g^{\epsilon})}_{\Omega_{x_k}}=g_k^{\epsilon} \qquad \mbox{and} \qquad \tilde{(f)}_{\Omega_{(x_k-\epsilon)}}=f_{k-1} \quad , \quad k=1,..,N
$$
Substituting this in eq.(26) we finally get an explicit expression for $F\cdot G^{\epsilon}$:
\begin{equation}
F\cdot G^{\epsilon}= fg^{\epsilon} + \sum_{k=1}^N \left[ g_k^{\epsilon}\Delta^{(F)}_{x_k}+ f_{k-1}\Delta^{(G)}_{x_k-\epsilon} \right]
\end{equation}
To obtain an explicit expression for $F\star G$ we still need to calculate the limit in eq.(21). Using eq.(29) we get for every $t \in \D(\bkR)$:
\begin{eqnarray}
\lim_{\epsilon \to 0^+} <F \cdot G^{\epsilon},t> &=& \lim_{\epsilon \to 0^+}\left\{ <f\cdot g^{\epsilon},t> +
\sum_{k=1}^N <g_k^{\epsilon}\Delta^{(F)}_{x_k},t>+ \sum_{k=1}^N <f_{k-1}\Delta^{(G)}_{x_k-\epsilon},t> \right\}\nonumber\\
&=&  \lim_{\epsilon \to 0^+}\int \, f(x)g(x+\epsilon) t(x) \, dx +\sum_{k=1}^N
 \lim_{\epsilon \to 0^+}<\Delta^{(F)}_{x_k},tg_k^{\epsilon}>  \nonumber \\
 && +\sum_{k=1}^N  \lim_{\epsilon \to 0^+}<\Delta^{(G)}_{x_k},t^{-\epsilon}f_{k-1}^{-\epsilon}> \\
&=& \int \, f(x)g(x) t(x) \, dx +\sum_{k=1}^N
<\Delta^{(F)}_{x_k},tg_k>  +\sum_{k=1}^N  <\Delta^{(G)}_{x_k},tf_{k-1}>
\nonumber \\
&=&<fg+\sum_{k=1}^N \left[ g_k \Delta^{(F)}_{x_k} +f_{k-1} \Delta^{(G)}_{x_k} \right] ,t> \nonumber
\end{eqnarray}
where we defined $t^{\epsilon}(x)=t(x+\epsilon)$, $f_k^{\epsilon}(x)=f_k(x+\epsilon)$ and used the fact that $t,g_k,f_k \in C^{\infty}(\bkR)$ and so they are uniformly continuous on an arbitrary compact interval.
We conclude from eqs.(21,30) that $F\star G \in {\mathcal D}'(\bkR)$ and is indeed given by eq.(24). Moreover, (24) is obviously of the form (7,8) and so $F\star G \in {\mathcal A}(\bkR)._{\Box}$\\

We proceed by studying some of the properties of the new product:
\\
\\
{\bf Theorem 3.16.} {\it The $\star$-product is (i) distributive, (ii) associative but (iii) non-commutative,
(iv) it reduces to the standard product of an infinitely smooth function by a distribution if either $F$ or $G$ belong to $C^{\infty}(\bkR)$ and (v) it reproduces the standard product of continuous functions.}\\
\\
{\bf Proof.}
\\
(i) Let us prove that the product is left distributive. Let $F,G,J \in {\mathcal A}(\bkR)$. We have:
\begin{equation}
(F+G)\star J= \lim_{\epsilon \to 0^+} (F+G) \cdot J^{\epsilon}
\end{equation}
We now prove that if $K\in {\mathcal A}(\bkR)$ and $V_{F+G}=V_F\cup V_G$ is such that $V_{F+G} \cap V_K= {\emptyset } $ then:
\begin{equation}
(F+G)\cdot K=F\cdot K+G\cdot K
\end{equation}
This identity is valid globally if it is valid locally.
Let $\{\Omega_w, w\in V_{F+G}\cup V_K\}$ be an open covering of $\bkR$ satisfying the conditions of Definition 3.4 for the $\cdot$ product of $F+G$ by $K$. Then eq.(32) is equivalent to:
\begin{equation}
[(F+G)\cdot K]_{\Omega_w}=[F\cdot K+G \cdot K]_{\Omega_w},\quad \forall w\in V_{F+G}\cup V_K
\end{equation}
which, in turn, follows from:
\begin{eqnarray}
[(F+G)\cdot K]_{\Omega_w} & = & (F+G)_w K_w=F_wK_w+G_wK_w \\
& = & [F\cdot K]_{\Omega_w}+[G\cdot K]_{\Omega_w}=[F\cdot K+G\cdot K]_{\Omega_w} \nonumber
\end{eqnarray}
where we used the definition of the $\cdot$ product (eq.(15)) and the fact that for each $w$ either $F_w,G_w\in C^{\infty}(\Omega_w)$ or $K_w\in C^{\infty}(\Omega_w)$. Making $K=J^{\epsilon}$ in eq.(32) and substituting in eq.(31), we get:
\begin{equation}
(F+G)\star J= \lim_{\epsilon \to 0^+} F \cdot J^{\epsilon} +\lim_{\epsilon \to 0^+} G \cdot J^{\epsilon}=F\star J +G\star J
\end{equation}
which proves that the product is left distributive.
Equivalently, one may prove that it is also right distributive.
\\
(ii) Let $F,G,J \in {\mathcal A}(\bkR)$. Then:
\begin{equation}
F= f+\sum_{w\in V_F} \Delta^{(F)}_{w} \quad ,\quad
G= g+\sum_{w\in V_G} \Delta^{(G)}_{w} \quad \mbox{and} \quad
J= j+\sum_{w\in V_J} \Delta^{(J)}_{w}
\end{equation}
where $f,g,j \in C_p^{\infty}(\bkR)$. To keep the notation simple let us redefine $V_F=V_G=V_J=$sing supp$\,F \,\cup \, $sing supp$\, G\, \cup \, $sing supp$\, J$. Let $N=\sharp V_F$ and let $V_F=\{x_1,...,x_N\}$ where $x_i<x_k$ for $i<k$. Then $F,G$ can be written as in eq.(22) with $f,g$ given by (23) and:
\begin{equation}
J= j+\sum_{k=1}^N \Delta^{(J)}_{x_k}
\end{equation}
(notice that $\Delta_{x_k}^{(F)}=0$ if $x_k \in V_F {\backslash }  \, $sing supp$\, F$ and equivalently for the other distributions) where $j$ can be written in terms of $C^{\infty}(\bkR)$-functions $j_k;\, k=0,..,N$, as:
\begin{equation}
j(x)  =  H(x_1-x)j_0(x) + \sum_{k=1}^{N-1} H(x-x_k) H(x_{k+1}-x) j_k(x) + H(x-x_N)j_N(x)
\end{equation}
A trivial calculation using eq.(24) shows that:
\begin{eqnarray}
 &&(F\star G)\star J  =  \left[ fg +\sum_{k=1}^N \left( g_k\Delta^{(F)}_{x_k}
+f_{k-1}\Delta^{(G)}_{x_k} \right) \right] \star J  \\
&=& fgj + \sum_{k=1}^N \left( g_kj_k\Delta^{(F)}_{x_k}
+f_{k-1}j_k\Delta^{(G)}_{x_k} +g_{k-1}f_{k-1}\Delta^{(J)}_{x_k} \right)
\nonumber
\end{eqnarray}
and also:
\begin{eqnarray}
 &&F\star (G\star J)  =  F\star \left[ gj +\sum_{k=1}^N \left( j_k\Delta^{(G)}_{x_k}+g_{k-1}\Delta^{(J)}_{x_k} \right) \right]  \\
&=& fgj + \sum_{k=1}^N \left( g_kj_k\Delta^{(F)}_{x_k}
+f_{k-1}j_k\Delta^{(G)}_{x_k} +g_{k-1}f_{k-1}\Delta^{(J)}_{x_k} \right)
\nonumber
\end{eqnarray}
Hence, the $\star$-product is associative.
\\
(iii)  Let $F,G$ be given by eq.(22,23). We have from eq.(24) that:
\begin{equation}
F\star G-G\star F= \sum_{k=1}^N \left[ \left( g_k -g_{k-1}\right) \Delta^{(F)}_{x_k}+\left( f_{k-1}-f_k \right) \Delta^{(G)}_{x_k} \right]
\end{equation}
and thus the $\star $-product is non-commutative. Notice nevertheless, that supp $(F\star G-G\star F) \subseteq V_F \cup V_G$, is a finite set.
\\
(iv) It follows directly from eq.(24) that if $F=f$ (or $G=f$) for some $f \in C^{\infty}(\bkR)$ then $F \star G=f G$ (respectively, $F \star G=f F$).
\\
(v) If $F,G \in ({\mathcal A}(\bkR) \cap C^0(\bkR))$ then, in the notation of eq.(22), $F=f$ and $G=g$. Hence, from eq.(24) $F\star G=fg$ thus
concluding the proof.
$_{\Box}$\\

Hence ${\mathcal A}(\bkR)$ is an associative algebra. To proceed we may endow ${\mathcal A}(\bkR)$ with the following bracket structure:
\begin{equation}
[F,G]=F\star G-G\star F
\end{equation}
the explicit form of the bracket being given by eq.(41).\\
\\
{\bf Theorem 3.17.} {\it The bracket (42) is a Lie bracket.}\\
\\
{\bf Proof.} The antisymmetric property of the bracket follows by construction while the linearity and the Jacobi identity are inherited from the left and right distributive and associative properties of the $\star$-product, respectively. It is also easy to check that the bracket satisfies the Leibniz rule with respect to the $\star$-product.$_{\Box}$ \\

Finally we prove that:\\
\\
{\bf Theorem 3.18.} {\it Let $\frac{d}{dx}$ be the usual distributional derivative in ${\mathcal D}'(\bkR)$. The restriction $\frac{d}{dx}: {\mathcal A}(\bkR) \to {\mathcal A}(\bkR)$ satisfies the Leibniz rule with respect to the $\star$-product.}\\
\\
{\bf Proof.} Let $F,G \in {\mathcal A}(\bkR)$. Then $F\star G \in {\mathcal A}(\bkR)$ and its derivative acts as ($F'=\frac{d}{dx} F$):
\begin{equation}
<(F\star G)',t>=-<F\star G,t'>=-\lim_{\epsilon \to 0^+} <F\cdot G^{\epsilon},t'> =\lim_{\epsilon \to 0^+} <(F\cdot G^{\epsilon})',t>
\end{equation}
Let us then prove that the product $\cdot$ satisfies the Leibniz rule. Let $F,J \in {\mathcal A}(\bkR)$ be such that $V_F\cap V_J= {\emptyset }
$. Let $\{\Omega_w,\, w\in V_F\cup V_J \}$ be an open covering of $\bkR$ satisfying the conditions of Definition 3.4. Locally we have:
\begin{eqnarray}
& & \left[ (F\cdot J)'\right]_{\Omega_w}=\left[ F'\cdot J+F\cdot J' \right]_{\Omega_w} \\
& \Longleftrightarrow &\left[ (F\cdot J)_{\Omega_w}\right]'=\left[ F'\cdot J \right]_{\Omega_w}+\left[F\cdot J' \right]_{\Omega_w}  \nonumber \\
&\Longleftrightarrow & (F_w J_w)' = F_w' J_w+F_w J_w' \nonumber
\end{eqnarray}
which is true because either $F_w$ or $J_w$ belongs to $C^{\infty}(\Omega_w)$. Hence the $\cdot$ product satisfies the Leibniz rule locally and
thus also globally.\\
Substituting this result in eq.(43) we get:
\begin{eqnarray}
<(F\star G)',t>& = & \lim_{\epsilon \to 0^+} <F\cdot \left(G^{\epsilon}\right)'+F'\cdot G^{\epsilon},t>  \\
& = & \lim_{\epsilon \to 0^+} <F\cdot \left(G'\right)^{\epsilon},t>+\lim_{\epsilon \to 0^+} <F'\cdot G^{\epsilon},t> \nonumber \\
&=& <F\star G',t>+<F'\star G,t> \nonumber
\end{eqnarray}
where we took into account that the translation and the derivative operators commute in $\D'(\bkR )$. Hence, $(F\star G)'=F\star G'+F'\star G$, which concludes the proof.$_{\Box}$\\
\\
{\bf Corollary 3.19.} The algebra $({\mathcal A}(\bkR),+,\star)$ is an associative (but non-commutative) differential algebra of generalized
functions.
\\

This concludes the study of the properties of the $\star$-product and of the associated algebra ${\mathcal A}(\bkR)$. As a simple example, let
us calculate the derivative of $H(x) \star H (x)$. We have:
$$
\frac{d}{dx} (H(x) \star H (x))= \delta(x) \star H(x) + H(x) \star \delta(x) = \delta(x) \star H (x) = \delta(x).
$$
which of course, is consistent with the fact that $H(x) \star H(x)= H(x)$. Notice that multiple products of the Heaviside step function are
sometimes used to exemplify a simpler version of the Schwartz impossibility result \cite{Colombeau3}. One typically finds that
differentiating the equation $H^n(x)=H(x)$ for different values of $n$ leads to contradictory results. This is not the case if $H^n(x)$ is
calculated using the product $\star$. The key point, as one can easily check, is that the $\star$ product is non-commutative.

To finish this section, let us point out that there are other associative but non-commutative products, related to the one introduced in
eq.(20), and also yielding a Lie bracket structure in ${\mathcal A}(\bkR)$. For instance:
\begin{eqnarray}
F\star_2 G & \equiv & \lim_{\epsilon \to 0^+} F\cdot G^{-\epsilon} \nonumber \\
F\star _3 G & \equiv &  \lim_{\epsilon \to 0^+} F^{\epsilon}\cdot G \\
F\star_4 G & \equiv & \lim_{\epsilon \to 0^+} F^{-\epsilon} \cdot G \nonumber
\end{eqnarray}
which are all related to the original product introduced in (20):
\begin{equation}
F\star_ 2G = F\star_3 G= G\star_4 F = G\star F
\end{equation}

Finally, we can also define a commutative product (or symmetric bracket) through the prescription:
\begin{equation}
F\star_5 G= \frac{1}{2} (F\star G+G\star F)
\end{equation}
which however is not associative.

\section{Linear differential equations with distributional coefficients}

In this section we consider the second problem presented in the Introduction. Let:
\begin{equation}
\sum_{i=0}^{n} a_i \psi^{(i)}  = f,
\end{equation}
be an ordinary, linear differential equation defined on $\bkR$. The coefficients are the complex valued functions $a_i \in C^{\infty}(\bkR)$ and $a_n (x) \not=0$, $\forall x \in \bkR$, and the non-homogeneous term is $f\in C^{\infty}(\bkR)$. As usual, $\psi^{(i)}$ denotes the $i$-th order derivative of $\psi$.
Let $\psi_U$ be the general solution of eq.(49). The subscript $U$ stands for {\it unconfined} (meaning the $\psi_U$ is the solution of (49) on the entire real line).
The generalized solution of eq.(49) (i.e. the solution in $\D'(\bkR)$) coincides with its classical solution (i.e. the one in $C^{\infty}(\bkR)$) \cite{Kanwal}. The designation $\psi_U$ will be used to denote both the distribution in $\D'(\bkR)$ as well as the associated infinitely smooth function.

Let us also consider an object of the form $\psi_C=H(x) \psi_U$ "confined" to the positive axis and reproducing the "unconfined" general solution for $x>0$.
Such a function is not a generalized (nor classical) solution of eq.(49). The aim of this section is to derive a new differential equation that displays the generalized solution (and only the generalized solution) $\psi_C$.

Notice that we are focusing on the particular case of Problem 2 (stated in the Introduction) where $\Omega=\bkR^+$. This is to keep the formulation as simple as possible. The extension to the general case where $\Omega =]a,b[$, $a<b \in \bkR \cup\{-\infty , +\infty \}$ is straightforward.

Let us then proceed. As a first step we substitute $\psi_C$ in eq.(49) and find the correction terms:
\begin{eqnarray}
\sum_{i=0}^{n} a_i \psi_C^{(i)}  & = & \sum_{i=0}^{n} a_i \left\{
\sum_{j=0}^i (^i_j) H^{(j)} \psi_U^{(i-j)} \right\}  \\
&=& \sum_{i=0}^{n} a_i H \psi_U^{(i)} +
\sum_{i=1}^{n} a_i
\sum_{j=1}^i (^i_j) H^{(j)} \psi_U^{(i-j)} \nonumber \\
&=& H f + \sum_{i \ge j =1}^{n} (^i_j) a_i  H^{(j)} \psi_U^{(i-j)} \nonumber
\end{eqnarray}
\\
{\bf Remark 4.1.} Equation (50) is not a closed equation for $\psi_C$, because it displays a distributional non-homogenous term which is dependent on the particular solution of eq.(49). For each solution of eq.(49) (associated to a complete set of boundary conditions $\psi_U^{(i)}(x_0), i=0,...,n-1$ given at some $x_0 \in \bkR^+$), this term can be re-written as a combination of Dirac deltas and their derivatives:
$$
\sum_{i \ge j =1}^{n} (^i_j) a_i(x)  H^{(j)}(x) \psi_U^{(i-j)}(x)= \sum_{i =1}^{n} \sum_{k=0}^{i-1} a_i(x)  \psi_U^{(k)} (0) \delta^{(i-1-k)}(x)
$$
and upon substitution in (50) we obtain a linear and ordinary differential equation, with a distributional non-homogenous term, that for a particular set of boundary conditions, displays a confined solution.
However, the general solution of this equation is still not confined to $\bkR^+$, but exhibits instead a jump discontinuity at $x=0$ \cite{Kanwal}. More details about ordinary differential equations with distributional non-homogeneous terms can be found in \cite{Estrada, Pan, Stakgold}. The related issue of distributional solutions of ordinary differential equations is studied in \cite{Littlejohn,Krall,Wiener}\\
\\
Our task here will be to re-write eq.(50) exclusively in terms of $\psi_C$ so that its general solution (and not only a particular solution) is of the form $H \psi_U$. With this in mind, we prove the following theorem:\\
\\
{\bf Theorem 4.2.} {\it If $\psi_U \in C^{\infty}(\bkR)$ and $\psi_C \in {\mathcal A}(\bkR)$ is such that $\psi_C=\psi_U$ on $\bkR^+$ then, in
the distributional sense:
\begin{equation}
\sum_{i\ge j=1}^{n} (^i_j) a_i H^{(j)} \psi_U^{(i-j)} =
\sum_{i\ge j=1}^{n} (^i_j) a_i H^{(j)} \star \psi_C^{(i-j)}.
\end{equation}}
\\
{\bf Proof.} The most general distribution $\psi_C \in {\mathcal A}(\bkR)$ satisfying $\psi_C=\psi_U$ on $\bkR^+$ is given by:
\begin{equation}
\psi_C=H \psi_U +F
\end{equation}
where $F \in {\mathcal A}(\bkR)$ is such that supp $F \subseteq \bkR^-\cup \{0\}$. Hence we have (for all $i \ge j=1,...,n$):
\begin{eqnarray}
<H^{(j)} \star {\psi}_C^{(i-j)} , t> &=& <H^{(j)} \star \left[H  {\psi}_U+F \right]^{(i-j)},t> \nonumber\\
&=&<H^{(j)} \star \left[H \star {\psi}_U \right]^{(i-j)},t>+<H^{(j)} \star F^{(i-j)},t> \nonumber \\
&=& <H^{(j)} \star {\psi}_U^{(i-j)},t>+\lim_{\epsilon \to 0^+}<H^{(j)} \cdot (F^{\epsilon})^{(i-j)},t> \nonumber \\
&=& <H^{(j)} {\psi}_U^{(i-j)} , t>
\end{eqnarray}
where in the last step we took into account that for $j\ge 1$, the sets supp $H^{(j)}=\{0\}$ and supp $(F^{\epsilon})^{(i-j)} \subseteq ]-\infty,-\epsilon]$ are disjoint.
The result (53) is valid for all $i\ge j=1,..,n$ and so eq.(51) holds.$_{\Box}$\\

We can now re-write eq.(50) exclusive in terms of $\psi_C$:
\begin{equation}
\sum_{i=0}^{n} a_i \psi_C^{(i)}  = H f+
\sum_{i\ge j=1}^{n} (^i_j) a_i H^{(j)} \star {\psi}_C^{(i-j)}
\end{equation}
which seems to be a suitable candidate for the differential equation we are looking for. That $\psi_C=H \psi_U$ is a solution of eq.(54) is obvious by construction. However, we still have to show that is its {\it only} solution:
\\
\\
{\bf Theorem 4.3.} {\it Let $\psi_U$ be the general solution of eq.(49). Then, the general solution of eq.(54) in ${\mathcal A}(\bkR)$ is given by $\psi_C=H \psi_U$.}\\
\\
{\bf Proof.} For a generic $\psi_C \in {\mathcal A}(\bkR)$ the support of the second term on the right hand side of eq.(54) is $\{0\}$. Hence,
on $\bkR^+$ and $\bkR^-$, eq.(54) reduces to eq.(49) and to the homogeneous equation associated to eq.(49), respectively. We conclude that any
solution of eq.(54) will be of the form:
\begin{equation}
\psi_C=\Delta +H _- \psi_{-} +H \psi_{+}
\end{equation}
where supp $\Delta \subseteq \{0\}$ (i.e $\Delta$ is a linear combination of the Dirac measure and its derivatives), $H_-(x)=H (-x)$ is the
reversed Heaviside step function $H$ and $\psi_{-},\psi_{+} \in C^{\infty}(\bkR)$ are the distributional (which coincide with the classical
\cite{Kanwal}) solutions of $\sum_{i=0}^n a_i \psi^{(i)}=0$ and of eq.(49), respectively. Substituting (55) in eq.(54) we get:
\begin{eqnarray}
&& \sum_{i=0}^{n} a_i\Delta ^{(i)} +
\sum_{i=0}^{n} \sum_{j=0}^{i} (^i_j) a_i \left( H_-^{(j)}  \psi_{-}^{(i-j)} +H^{(j)} \psi_{+}^{(i-j)} \right)  \\
&& = H f + \sum_{i \ge j =1}^n (^i_j) a_i  H^{(j)} \psi_{+}^{(i-j)} \nonumber \\
&\Longleftrightarrow &
\sum_{i=0}^{n} a_i\Delta ^{(i)} +
\sum_{i=0}^{n} a_i \left( H_- \psi_{-}^{(i)}(x) +H \psi_{+}^{(i)} \right)  \nonumber \\
&& +\sum_{i \ge j=1}^{n} (^i_j) a_i \left( H^{(j)}  \psi_{+}^{(i-j)} -H^{(j)} \psi_{-}^{(i-j)} \right) = H f + \sum_{i \ge j =1}^n (^i_j) a_i
H^{(j)}  \psi_{+}^{(i-j)} \nonumber
\end{eqnarray}
Separating the terms that involve the delta distribution or its derivatives from those that do not, we get:
\begin{eqnarray}
&& \sum_{i=0}^{n} a_i \left( H_- \psi_{-}^{(i)} +H \psi_{+}^{(i)} \right) =H f \\
& \Longleftrightarrow & \qquad \left\{ \begin{array}{lll}
\sum_{i=0}^{n} a_i(x)  \psi_{+}^{(i)}(x) & = & f(x)   , \quad x > 0 \\
\\
\sum_{i=0}^{n} a_i(x) \psi_{-}^{(i)}(x) & = & 0  \quad  , \quad x<0
\end{array} \right. \nonumber
\end{eqnarray}
(confirming what we already knew about $\psi_-,\psi_+$) and:
\begin{equation}
\sum_{i=0}^{n} a_i\Delta ^{(i)} -
\sum_{i \ge j=1}^{n} (^i_j) a_i  H^{(j)} \psi_{-}^{(i-j)} =0
\end{equation}
In eq.(58) the term proportional to the highest order derivative of $\delta$ is given by
$ a_n \Delta^{(n)}$ (if not zero, $a_n \Delta^{(n)}$ is proportional to at least $\delta^{(n)}$). This term cannot be cancelled by a combination of lower order derivatives of $\delta$. Since $a_n \not=0$
we get $\Delta^{(n)} =0$ and so $\Delta=0$. Hence, eq.(58) reduces to:
\begin{equation}
\sum_{i \ge j=1}^{n} (^i_j) a_i  H^{(j)} \psi_{-}^{(i-j)} =0 \Longleftrightarrow \sum_{j=1}^n \sum_{i=j}^{n} (^i_j) a_i  H_-^{(j)} \psi_{-}^{(i-j)} =0
\end{equation}
From eq.(57) we find that $\sum_{i=0}^n a_i H_-\psi_-^{(i)}=0 $, $\forall x\in \bkR$ and by adding this last equation to eq.(59) we finally get:
\begin{equation}
\sum_{j=0}^n \sum_{i=j}^{n} (^i_j) a_i  H_-^{(j)} \psi_{-}^{(i-j)} =0 \Longleftrightarrow \sum_{i=0}^n a_i \frac{d^i}{dx^i} \left(
H_-\psi_{-} \right)=0 \Longleftrightarrow
\sum_{i=0}^n a_i \phi^{(i)}=0
\end{equation}
where $\phi=H_- \psi_- $ satisfies (by construction) the boundary conditions $\phi^{(i)}(z_0)=0$, $i=0,...n-1$ at an arbitrary $z_0>0$. For these boundary conditions the unique solution is $\phi =0$, which implies that $\psi_-(x<0)=0$.

Assembling all these results, we finally get: $\psi_C =H \psi_+$, where $\psi_+$ satisfies eq.(57) (or equivalently eq.(49)), which proves the theorem.$_{\Box}$
\\
\\
We conclude that eq.(54) provides a global formulation for the differential problem described by eq.(49) and confined to the positive axis, i.e. it yields a solution of Problem 2 for the case where $\Omega=\bkR^+$. The extension to the more general case where $\Omega$ is an arbitrary interval is straightforward.

\section*{Acknowledgement}
We would like to thank Carlos Sarrico for useful comments and insights into the subject. Research was supported by the grants POCTI/0208/2003
and PTDC/MAT/69635/2006 of the Portuguese Science Foundation.

\end{document}